\documentclass[12pt]{amsart}
\usepackage{amsmath,amssymb,amsthm,amscd}
\usepackage{graphicx}
\usepackage{color}
\usepackage{enumitem}
\usepackage[plainpages=false]{hyperref}

\newcommand{\beq}{\begin{equation}}
\newcommand{\eeq}{\end{equation}}
\newcommand{\beqs}{\begin{eqnarray*}}
	\newcommand{\eeqs}{\end{eqnarray*}}
\newcommand{\beqn}{\begin{eqnarray}}
\newcommand{\eeqn}{\end{eqnarray}}
\newcommand{\beqa}{\begin{array}}
	\newcommand{\eeqa}{\end{array}}

\newcommand{\re}{\text{Re}}

\DeclareMathOperator{\sym}{Sym}
\DeclareMathOperator{\herm}{Herm}
\DeclareMathOperator{\tr}{tr}

\newtheorem{prop}{Proposition}[section]
\newtheorem{thm}[prop]{Theorem}

\newtheorem{cor}[prop]{Corollary}
\newtheorem{rem}[prop]{Remark}

\begin{document}
 
\title[Bernstein Theorem for the Self-Shrinking $J$-Equation]
{A Bernstein Theorem for the Self-Shrinking $J$-Equation and Some Generalizations}

\author{Yiyang Pan}
\address [Yiyang Pan]{School of Mathematical Sciences and LPMC, Nankai University, Tianjin, 300071, P. R. China}
\email{yiyang.pan.math@foxmail.com}
\author{Wenlong Wang}
\address [Wenlong Wang]{School of Mathematical Sciences and LPMC, Nankai University, Tianjin, 300071, P. R. China}
\email{wangwl@nankai.edu.cn}

\begin{abstract}
We prove that every entire smooth strictly plurisubharmonic solution of the self-shrinking $J$-equation on $\mathbb{C}^n$ is a quadratic polynomial. This removes the asymptotic lower bound assumption on the complex Hessian in \cite[Theorem 4]{HJ}. The result also recovers the corresponding real rigidity theorem in \cite[Theorem 1.1]{HOW} as a special case. More generally, our method applies to self-shrinking equations associated with a broad class of fully nonlinear elliptic operators satisfying suitable structural conditions, including the adjacent inverse complex Hessian quotient operators $-\sigma_{k-1}/\sigma_{k}$, $1\leq k\leq n$. For the remaining quotients $-\sigma_l/\sigma_k$ with $l<k-1$, we prove the same rigidity conclusion under the additional assumption that the solution is proper.

\vspace{0.15cm}

\vspace{0.15cm}
\end{abstract}
\maketitle

\section{Introduction}
The $J$-flow was introduced by Donaldson \cite{Donaldson} in the moment-map framework, and independently by Chen \cite{CXX} as the gradient flow of the $J$-functional appearing
in his formula for the Mabuchi energy. Its stationary points satisfy the $J$-equation, an important fully nonlinear geometric PDE closely related to the study of canonical metrics. Over the past two decades, a substantial theory has been developed in the compact K\"{a}hler setting; see, for example, \cite{CXX2,Weinkove1,Weinkove2,SW,FLM,LS,CS,Sz,CG,DP,Song} and the references therein. In this paper, we study a self-similar version of the $J$-equation on $\mathbb{C}^n$ and prove a Bernstein theorem for its entire solutions.

We use the following conventions throughout the paper. On $\mathbb R^m$, we write $x=(x^1,\ldots,x^m)$ for the standard coordinates; $D$ denotes the Euclidean gradient, and ``$\cdot$'' denotes the Euclidean inner product. On $\mathbb C^n$, we use the standard complex coordinates $z=(z^1,\ldots,z^n)$ and identify $\mathbb C^n$ with $\mathbb R^{2n}$ by writing $z^i=x^i+\sqrt{-1}\,x^{n+i}$ for $i=1,\ldots,n$. We can now state our main result for the self-shrinking $J$-equation.

\begin{thm}\label{Jthm}
Let $u:\mathbb{C}^n\to\mathbb{R}$ be a smooth strictly plurisubharmonic function satisfying the self-shrinking $J$-equation
\begin{equation}\label{Jeq}
-\sum_{i=1}^n\frac{1}{\lambda_i}=\frac{1}{2}x\cdot Du-u,
\end{equation}
where $\lambda_1,\ldots,\lambda_n$ are the eigenvalues of the complex Hessian $\partial\bar\partial u$. Then $u$ is a quadratic polynomial. 
\end{thm}
Any solution $u$ of \eqref{Jeq} gives rise to the self-shrinking solution \[v(x,t)=-tu\left(\frac{x}{\sqrt{-t}}\right)\] of the parabolic complex inverse trace equation
\begin{equation}\label{Jflow}
\partial_tv=-\sum^n_{i=1}\frac{1}{\lambda_i}
\end{equation}
on $\mathbb{C}^{n}\times(-\infty,0)$, where now $\lambda_1,\ldots,\lambda_n$ are the eigenvalues of $\partial\bar\partial v$. In terms of K\"{a}hler forms, this equation can be written as
\[
\partial_{t}v=-\tr_\omega\omega_0,
\]
where $\omega=\sqrt{-1}\partial\bar\partial v$ and $\omega_0=\sqrt{-1}\partial\bar\partial |z|^2$. Up to the normalizing constant in the compact setting, \eqref{Jflow} is the natural $\mathbb C^n$ version of the $J$-flow with the fixed Euclidean background form $\omega_0$.

Theorem \ref{Jthm} removes the additional asymptotic lower bound assumption in Han--Jin
\cite[Theorem 4]{HJ}. There, the same conclusion was proved under the assumption that, for some fixed $\delta>0$,
\[
\partial\bar\partial u(x)\geq\frac{\sqrt{2n-1+\delta}}{|x|}I_n
\]
holds for sufficiently large $|x|$. 

Theorem \ref{Jthm} also recovers the corresponding real rigidity theorem of Huang, Ou, and the second author \cite[Theorem 1.1]{HOW} as a special case. After an elementary normalization of the potential, their theorem asserts that every entire smooth strictly convex solution $w$ on $\mathbb R^n$ of
\begin{equation}\label{RJeq}
-\sum_{i=1}^n \frac{1}{\mu_i}=\frac{1}{2}x\cdot Dw-w
\end{equation}
is quadratic, where $\mu_1,\ldots,\mu_n$ are the eigenvalues of $D^2w$. To see the reduction, set $u(z)=2w(\re z)$. Then $u$ is strictly plurisubharmonic on $\mathbb C^n$ and solves \eqref{Jeq}, so Theorem \ref{Jthm} yields the quadraticity of $w$. The geometric interpretation of \eqref{RJeq} is explained in Section 2 of \cite{HOW}: for a smooth strictly convex solution $w$ on $\mathbb R^n$, the gradient graph
\[
\Sigma=\{(x,Dw(x))\mid x\in\mathbb{R}^n\}
\]
is a Lagrangian self-shrinker in $\mathbb{R}^{2n}$ equipped with the degenerate metric 
\[
g=\sum_{i=1}^n (dx^{n+i})^2.
\]

In fact, the proof of Theorem \ref{Jthm} does not rely on the specific form of the $J$-operator. Rather, it applies to a broader class of fully nonlinear elliptic operators satisfying suitable structural conditions.

To state the general result, let $\sym(n)$ denote the space of real symmetric $n\times n$ matrices, and let $F$ be a $C^1$ function defined on an open subset $U\subset\sym(n)$. For $B\in U$, let $DF(B)$ denote the linearized coefficient matrix of $F$ at $B$, defined by
\begin{equation*}
\frac{d}{dt}\bigg|_{t=0}F(B+tK)=\tr\left(DF(B)K\right)\qquad\text{for every $K\in\sym(n)$.}
\end{equation*}
We assume that $F$ satisfies the following conditions:
\begin{itemize}[leftmargin=0.8cm]
\item[(i)] \emph{Ellipticity}: $DF(B)$ is positive definite for every $B\in U$;
\item[(ii)] \emph{Upper boundedness}: there exists a constant $c\in\mathbb{R}$ such that
\[
F(B)\le c
\qquad\text{for every }B\in U;
\]
\item[(iii)] \emph{Coefficient control}: there exists a continuous nonnegative function $G$ on $(-\infty,c]$ such that
\[
\tr (DF(B))\le G(F(B))
\qquad\text{for every }B\in U.
\]
\end{itemize}

We obtain the following general rigidity theorem.
\begin{thm}\label{Gthm}
Assume that $F$ satisfies conditions \textup{(i)}--\textup{(iii)}.
Let $u\in C^3(\mathbb R^n)$ satisfy
$D^2u(\mathbb R^n)\subset U$ and
\begin{equation}\label{Geq}
F(D^2u)=\frac12 x\cdot Du-u.
\end{equation}
Then $u$ is a quadratic polynomial.
\end{thm}

Although the $J$-operator is naturally expressed in terms of the complex Hessian, it can also be regarded as an operator on real symmetric matrices and satisfies the structural conditions of Theorem \ref{Gthm}. The same is true for the inverse complex Hessian quotient operators $-\sigma_{k-1}/\sigma_k$, which include the $J$-operator as the case $k=n$. We recall the relevant notation.

For $1\leq k\leq n$, let $\sigma_k$ denote the $k$-th elementary symmetric polynomial,
\[
\sigma_k(\lambda)
=
\sum_{i_1<\cdots<i_k}
\lambda_{i_1}\cdots\lambda_{i_k},
\qquad
\lambda=(\lambda_1,\ldots,\lambda_n)\in\mathbb{R}^n,
\]
and use the convention $\sigma_0=1$. Set
\[
\Gamma_k
=
\{\lambda\in\mathbb{R}^n\mid
\sigma_j(\lambda)>0,\,1\leq j\leq k\}.
\]

For a real symmetric or Hermitian matrix $A$, let $\lambda(A)$ denote its eigenvalues. A $C^2$ function $u$ on $\Omega\subset\mathbb{R}^n$ is called $k$-convex if $\lambda(D^2u)\in\Gamma_k$ in $\Omega$. Similarly, a $C^2$ function $u$ on $\Omega\subset\mathbb{C}^n$ is called $k$-plurisubharmonic if $\lambda(\partial\bar\partial u)\in\Gamma_k$ in $\Omega$.

For $0\leq l<k$, define the real and complex inverse Hessian quotient operators by
\[
-\frac{\sigma_l}{\sigma_k}(D^2u)
=
-\frac{\sigma_l(\lambda(D^2u))}{\sigma_k(\lambda(D^2u))},
\qquad
-\frac{\sigma_l}{\sigma_k}(\partial\bar\partial u)
=
-\frac{\sigma_l(\lambda(\partial\bar\partial u))}{\sigma_k(\lambda(\partial\bar\partial u))},
\]
which are elliptic on the $k$-convex and $k$-plurisubharmonic classes, respectively.

As an application of Theorem \ref{Gthm}, we obtain the following Bernstein theorem.
\begin{thm}\label{CIHQthm}
Let $1\leq k\leq n$, and let $u:\mathbb{C}^n\to\mathbb{R}$ be a smooth $k$-plurisubharmonic function satisfying
\begin{equation*}
-\frac{\sigma_{k-1}}{\sigma_k}(\partial\bar\partial u)
=\frac{1}{2}x\cdot Du-u.
\end{equation*}
Then $u$ is a quadratic polynomial.
\end{thm}
The corresponding real statement follows immediately.

\begin{cor}\label{RIHQthm}
Let $1\leq k\leq n$, and let $u:\mathbb{R}^n\to\mathbb{R}$ be a smooth $k$-convex function satisfying
\begin{equation*}
-\frac{\sigma_{k-1}}{\sigma_k}(D^2u)
=\frac{1}{2}x\cdot Du-u.
\end{equation*}
Then $u$ is a quadratic polynomial.
\end{cor}
\begin{rem}
If condition \textup{(ii)} in Theorem \ref{Gthm} is removed, the conclusion can fail. For instance, the linear equation
\begin{equation*}
\Delta u=\frac{1}{2}x\cdot Du-u
\end{equation*}
satisfies condition \textup{(iii)} with $G\equiv n$, but admits a non-quadratic entire smooth solution; see Remark 2.1 in \cite{WWL1} for the explicit form. 
\end{rem}

On the other hand, the operators $-\sigma_{l}/\sigma_{k}$ with $k\geq l+2$ satisfy condition \textup{(ii)} but not condition \textup{(iii)}. Thus Theorem \ref{Gthm} does not apply to these operators. Nevertheless, under the additional assumption that the solution is proper, the same rigidity conclusion holds for the full family $-\sigma_l/\sigma_k$, $0\leq l<k\leq n$. 
\begin{thm}\label{GCIHQthm}
Let $0\leq l<k\leq n$, and let $u:\mathbb{C}^n\to\mathbb{R}$ be a smooth $k$-plurisubharmonic function satisfying
\begin{equation*}
-\frac{\sigma_{l}}{\sigma_k}(\partial\bar\partial u)
=\frac{1}{2}x\cdot Du-u.
\end{equation*}
Suppose in addition that $u$ is proper. Then $u$ is a quadratic polynomial.
\end{thm}
The same proof shows that every smooth proper \(k\)-convex solution of the corresponding real Hessian quotient equation on \(\mathbb R^n\) is a quadratic polynomial.

Related self-similar equations, with other fully nonlinear elliptic operators in place of the $J$-operator, also arise naturally in geometry and have been well studied.

When the operator is $\sum^n_{i=1}\arctan\lambda_i(D^2u)$, one obtains the potential equation for Lagrangian self-shrinkers in $\mathbb R^{2n}$ \cite{CCH}. This rigidity problem was studied in \cite{CCH,CCH2,HW}, and the corresponding Bernstein theorem was established by Chau--Chen--Yuan \cite{CCY}. Its complex counterpart, with operator $\sum^n_{i=1}\arctan\lambda_i(\partial\bar\partial u)$, is the potential equation on $\mathbb C^n$ for self-shrinking solutions to the line bundle mean curvature flow introduced by Jacob--Yau \cite{JY}. The corresponding Bernstein theorem was proved by Han--Jin \cite{HJ}. We also note that this case falls within the scope of Theorem \ref{Gthm}.

When the operator is $\log\det D^2u$, one obtains the potential equation for Lagrangian self-shrinkers in the pseudo-Euclidean space $(\mathbb R^{2n},g)$, where $g=\sum^n_{i=1}dx^idx^{n+i}$ \cite{H}. This rigidity problem was studied in \cite{HW,CCY}, and the Bernstein theorem was established by Ding--Xin \cite{DX}. Its complex counterpart, with operator $\log\det\partial\bar\partial u$, is the potential equation for shrinking K\"{a}hler--Ricci solitons on $\mathbb C^n$ \cite{CCY}. Under the completeness assumption, rigidity was proved by Drugan--Lu--Yuan \cite{DLY}. Without this assumption, the Bernstein problem was settled in complex dimension one by the second author \cite{WWL1}. In higher dimensions, the general Bernstein problem remains open, although rigidity has been established for $S^1$-symmetric solutions and for convex proper solutions \cite{WWL2}.

For the operator $\log\frac{\sigma_k}{\sigma_l}(D^2u)$, with $k>l$, and more generally for a class of operators satisfying certain structural conditions, rigidity for strictly convex solutions was proved by the second author \cite{WWL3}. These conditions, however, are not well suited to the present setting.

The proof of Theorem \ref{Gthm} builds on the barrier argument of Chau--Chen--Yuan \cite{CCY} for the self-similar term $\phi=\frac{1}{2}x\cdot Du-u$. They observed that, for each of the three equations considered in their work, $\phi$ satisfies a second-order elliptic equation without a zeroth-order term. Using this observation, they constructed a barrier function to force the global supremum of $\phi$ to be attained at some point. The strong maximum principle then implies that $\phi$ is constant, from which the quadratic conclusion follows. This strategy was also used in \cite{WWL3}, \cite{HOW} and \cite{HJ}. 

The main difficulty in the present setting lies in the barrier step, namely, forcing the self-similar term $\phi$ to attain its global maximum. In \cite{CCY,HJ}, the coefficient control required for this step is obtained directly from additional assumptions or from properties specific to those operators. For the self-shrinking $J$-equation, we observe instead that the available coefficient control is tied to $\phi$ itself; this feature is formulated as condition (iii) in Theorem \ref{Gthm}. We therefore construct a family of auxiliary functions on expanding balls, built around $\phi$. 
The comparison with the value at the origin, together with the upper boundedness condition (ii), forces the values of $\phi$ at the maximizers of these auxiliary functions to remain in a fixed interval, so that condition (iii) provides the needed coefficient control at those points. The resulting maximum-principle estimates then keep these points in a fixed compact set as the balls exhaust $\mathbb R^n$, thereby completing the barrier step. 

For the non-adjacent quotients $-\sigma_l/\sigma_k$ with $k\geq l+2$, the coefficient-control mechanism above is no longer available. When $u$ is proper, however, we use the solution itself as a penalizing function to force $\phi$ to attain its global maximum. A difficulty in carrying out this argument is that $u$ is not known a priori to be a supersolution of the associated linearized equation.

{\it Organization.} In Section \ref{proofofmainthm}, we prove Theorem \ref{Gthm}. In Section \ref{verification}, we verify that the inverse complex Hessian quotient operators $-\sigma_{k-1}/\sigma_{k}$ on the $k$-plurisubharmonic class satisfy the hypotheses of Theorem \ref{Gthm}, thereby proving Theorem \ref{CIHQthm}. We also show that the deformed Hermitian--Yang--Mills phase operator and several further examples fall within the scope of Theorem \ref{Gthm}. In Section \ref{proofofginverseqthm}, we prove Theorem \ref{GCIHQthm}.

\section{Proof of Theorem \ref{Gthm}}\label{proofofmainthm}
The proof is divided into three steps. We first derive an elliptic equation for the self-similar term. 
\subsection{Equation for the self-similar term}\label{subsecEquation for the self-similar term}
Set
\begin{equation}\label{phase1}
\phi(x):=\frac{1}{2}x\cdot Du(x)-u(x).
\end{equation}
By \eqref{Geq}, we also have
\begin{equation}\label{phase2}
\phi=F(D^2u).
\end{equation}
We write $F^{ij}$ for the $(i,j)$-entry of $DF(D^2u)$. 
Differentiating \eqref{phase2} in the $x_s$-direction gives
\begin{equation}\label{diffu}
F^{ij}u_{ijs}=\phi_{s}.
\end{equation}
On the other hand, differentiating \eqref{phase1} twice yields
\begin{equation}\label{diffphase}
\phi_{ij}=\frac{1}{2}x\cdot Du_{ij}.
\end{equation}
Combining  \eqref{diffu} and \eqref{diffphase}, we obtain
\begin{equation}\label{Lphase}
F^{ij}\phi_{ij}-\frac{1}{2}x\cdot D\phi=0.
\end{equation}
Thus $\phi$ satisfies an elliptic equation without a zeroth-order term. In what follows, we write
\begin{equation*}
\mathcal{L}=F^{ij}\partial_{ij}-\frac{1}{2}x\cdot D.
\end{equation*}

Applying $\mathcal{L}$ to $|x|^2$, we get
\begin{equation*}
\mathcal{L}|x|^2=2\tr\left(DF(D^2u)\right)-|x|^2.
\end{equation*}
By the coefficient control condition (iii),
\begin{equation}\label{tracecontrol}
\tr \left(DF(D^2u)\right)\leq G(F(D^2u))=G(\phi).
\end{equation}
It follows that
\begin{equation}\label{Lx2}
\mathcal{L}|x|^2\leq 2G(\phi)-|x|^2.
\end{equation}

\subsection{Barrier construction}

For $R>0$, let $B_{R}$ be the open ball centered at the origin with radius $R$. On $\overline{B}_R$, consider the auxiliary function 
\[
f(x)=e^{\phi(x)}\left(R^2-|x|^2\right),
\]
and set 
\[
M(R)=\max_{\overline{B}_R} f(x).
\]
We first show that every maximizer of $f$ lies in a fixed compact
set independent of $R$. This localization will then yield the
required upper bound for $M(R)$. 

Since $f=0$ on $\partial B_R$, the maximum is attained at some point $y\in B_R$. At this point,
\[
e^{\phi(y)}\left(R^2-|y|^2\right)\geq e^{\phi(0)}R^2.
\]
It follows that $\phi(y)\geq\phi(0)$. Since $\phi=F(D^2u)$, condition (ii) implies $\phi(y)\leq c$. Hence,
\begin{equation}\label{leveldis}
R^2-|y|^2\geq e^{\phi(0)-\phi(y)}R^2\geq e^{\phi(0)-c}R^2.
\end{equation}
At $y$, write $A=(a^{ij})=DF(D^2u(y))$. Since $y$ is an interior maximizer of $f$, we have 
\begin{equation*}
Df(y)=0,\qquad a^{ij}f_{ij}(y)\leq 0.
\end{equation*}
Since $Df(y)=0$, the drift term in $\mathcal Lf$ vanishes at $y$. Hence
\begin{equation}\label{maximumL}
\mathcal Lf(y)\leq 0.
\end{equation}
Expanding $Df(y)=0$, we obtain
\begin{equation}\label{gradienteq}
\left(R^2-|y|^2\right)D\phi(y)=2y.
\end{equation}
Expanding $\mathcal Lf(y)$, we get
\begin{equation}\label{Lfexpand}
\begin{split}
\mathcal Lf(y)&=\left(R^2-|y|^2\right)\mathcal{L}(e^{\phi})(y)+e^{\phi(y)}\mathcal L\left(R^2-|x|^2\right)(y)-4e^{\phi(y)}a^{ij}\phi_i(y)y_j\\
&=:\mathrm{I}+\mathrm{II}+\mathrm{III}.
\end{split}
\end{equation}
By \eqref{Lphase},
\begin{equation}\label{Iestimate}
\begin{split}
\mathrm{I}&=\left(R^2-|y|^2\right)e^{\phi(y)}\left(\mathcal{L}\phi(y)+a^{ij}\phi_i(y)\phi_j(y)\right)\\
&=\left(R^2-|y|^2\right)e^{\phi(y)}a^{ij}\phi_i(y)\phi_j(y)\geq 0.
\end{split}
\end{equation}
By \eqref{Lx2},
\begin{equation}\label{IIestimate}
\mathrm{II}\geq e^{\phi(y)}\left(|y|^2-2G(\phi(y))\right).
\end{equation}
Using \eqref{gradienteq}, we rewrite $\mathrm{III}$ as
\begin{equation}\label{expressionofIII}
\mathrm{III}=-\frac{8e^{\phi(y)}}{R^2-|y|^2}a^{ij}y_iy_j.
\end{equation}
By the positivity of $A$ and the trace estimate \eqref{tracecontrol}, we have
\begin{equation}\label{quadraticcontrol}
a^{ij}y_iy_j\leq (\tr A)|y|^2\leq G(\phi(y))|y|^2\leq G(\phi(y))R^2.
\end{equation}
Using \eqref{leveldis} and \eqref{quadraticcontrol} in \eqref{expressionofIII}, we get
\begin{equation}\label{IIIestimate}
\mathrm{III}\geq-8e^{\phi(y)+c-\phi(0)}G(\phi(y)).
\end{equation}
Substituting \eqref{Iestimate}, \eqref{IIestimate}, and \eqref{IIIestimate} into \eqref{Lfexpand} and using \eqref{maximumL}, we obtain
\begin{equation*}
0\geq\mathcal Lf(y)\geq  e^{\phi(y)}\left[|y|^2-\left(2+8e^{c-\phi(0)}\right)G(\phi(y))\right].
\end{equation*}
Therefore, 
\begin{equation}\label{yestimate}
|y|^2\leq \left(2+8e^{c-\phi(0)}\right)G(\phi(y)).
\end{equation}
By continuity of $G$ on $[\phi(0),c]$, let
\begin{equation*}
G_0:=\max_{s\in[\phi(0),c]}G(s),\qquad R_0:=\sqrt{(2+8e^{c-\phi(0)})G_0}.
\end{equation*}
Since $\phi(y)\in [\phi(0),c]$, we have $G(\phi(y))\leq G_0$. Then \eqref{yestimate} gives 
\[
|y|\leq R_0.
\]
Therefore, $\phi(y)\leq \max_{\overline B_{R_0}}\phi$. It follows that
\begin{equation}\label{keyestimate}
M(R)\leq e^{\max_{\overline B_{R_0}}\phi}R^2.
\end{equation}
Here $R_0$ is independent of $R$, and the estimate holds for every $R>0$.

\subsection{Constancy of $\phi$ and quadraticity of $u$}
We show that estimate \eqref{keyestimate} forces 
\[
\phi\leq \max_{\overline B_{R_0}}\phi\quad\text{on}\ \mathbb R^n.
\] 
Otherwise, there exists $p\in\mathbb R^n$ such that $\phi(p)>\max_{\overline B_{R_0}}\phi$. Then, for every $R>|p|$, 
\begin{equation*}
M(R)\geq e^{\phi(p)}\left(R^2-|p|^2\right).
\end{equation*}
By comparing the leading $R^2$-terms, this contradicts \eqref{keyestimate} for sufficiently large $R$. Hence $\phi$ attains its global maximum in $\overline B_{R_0}$.  Applying the strong maximum principle to \eqref{Lphase}, we conclude that $\phi$ is constant.  

As $\phi$ is constant, \eqref{diffphase} gives
\[
x\cdot Du_{ij}=0.
\]
Thus, for each fixed $\theta\in\partial B_1$, the function $r\mapsto u_{ij}(r\theta)$ is constant for $r>0$. Since $u_{ij}$ is continuous, letting $r\to 0$ gives $u_{ij}(x)\equiv u_{ij}(0)$. Consequently, $D^2u$ is constant on $\mathbb R^n$, and $u$ is a quadratic polynomial.

\section{Verification for $-\sigma_{k-1}/\sigma_{k}$ and other operators}\label{verification}

We now verify the hypotheses of Theorem \ref{Gthm} for $-\sigma_{k-1}/\sigma_{k}$ and several other operators. Since some of these operators are formulated in the complex setting, we first explain how to regard them as operators on real symmetric matrices.

\subsection{The real--complex correspondence}
Let $\herm(n)$ denote the space of $n\times n$ Hermitian matrices.
Define the linear map 
\[
\mathcal H:\sym(2n)\rightarrow\herm(n)
\]
as follows. For \(M=(M_{ij})\in\sym(2n)\), set
\begin{equation*}
\bigl(\mathcal H(M)\bigr)_{i\bar j}
=
\frac14\Bigl(M_{ij}+M_{n+i,n+j}
+\sqrt{-1}\bigl(M_{i,n+j}-M_{n+i,j}\bigr)\Bigr),
\quad\  1\le i,\,j\le n.
\end{equation*}
Equivalently, if
\[
M=
\begin{pmatrix}
A & B\\
B^{\mathrm T} & C
\end{pmatrix},
\qquad
A,\,C\in \sym(n),\ B\in M_n(\mathbb R),
\]
then
\[
\mathcal H(M)=\frac14\left(A+C+\sqrt{-1}\left(B-B^{\mathrm T}\right)\right).
\]
Under the identification of \(\mathbb C^n\) with \(\mathbb R^{2n}\) given by
\(z^i=x^i+\sqrt{-1}\,x^{n+i}\), \(i=1,\ldots,n\), one has
\[
\mathcal H(D^2u)=\partial\bar\partial u
\]
for every \(C^2\) function \(u\) on \(\mathbb C^n\).

Let \(\Phi\) be a \(C^1\) function defined on an open subset \(V\subset\herm(n)\). Then \(\Phi\) induces a $C^1$ function $\widetilde{\Phi}$ on $\widetilde V:=\mathcal H^{-1}(V)\subset\sym(2n)$ by 
\[
\widetilde{\Phi}(M)=\Phi(\mathcal H(M)),\qquad M\in\widetilde V.
\]
For $Q\in V$, let $D\Phi(Q)\in\herm(n)$ be the unique matrix satisfying
\[
\frac{d}{dt}\bigg|_{t=0}\Phi\left(Q+tP\right)=\tr\left(D\Phi(Q)P\right)\qquad\text{for every $P\in \herm(n)$}.
\]
Since $\mathcal H$ is linear, the chain rule gives, for $M\in\widetilde V$ and $K\in\sym(2n)$, 
\begin{equation*}
\tr\big(D\widetilde\Phi\left(M\right)K\big)=\tr\big(D\Phi(\mathcal H(M))\mathcal H(K)\big).
\end{equation*}
Writing
\begin{equation*}
D\Phi(\mathcal H(M))=X+\sqrt{-1}Y,
\end{equation*} 
where $X$ is real symmetric and $Y$ is real skew-symmetric, we obtain 
\begin{equation*}
D\widetilde\Phi\left(M\right)=\frac{1}{4}\begin{pmatrix}
X & Y\\
-Y & X
\end{pmatrix}.
\end{equation*}
Consequently,
$D\widetilde\Phi\left(M\right)$ is positive definite if and only if $D\Phi(\mathcal H(M))$ is positive definite. Moreover, 
\begin{equation}\label{tracerelation}
\tr\big(D\widetilde\Phi\left(M\right)\big)=\frac{1}{2}\tr\big(D\Phi(\mathcal H(M))\big).
\end{equation} 

\subsection{Inverse Hessian quotient operator}\label{IHQP} For a Hermitian matrix $Q$ with $\lambda(Q)\in\Gamma_k$, set 
\[
\Phi_k(Q) := -\frac{\sigma_{k-1}(\lambda(Q))}{\sigma_k(\lambda(Q))}.
\]
By the real--complex correspondence, the induced real operator $\widetilde\Phi_k$ satisfies 
\[
\widetilde\Phi_k(D^2u)=\Phi_k(\mathcal H(D^2u))=-\frac{\sigma_{k-1}}{\sigma_k}(\partial\bar\partial u).
\] 
We now verify that $\widetilde\Phi_k$ satisfies the structural conditions of Theorem \ref{Gthm}. 

Since $\Phi_k<0$ on $\Gamma_k$, the upper boundedness condition \textup{(ii)} holds with $c=0$. We now compute the linearized coefficients. By unitary invariance, it suffices to compute at a diagonal matrix \(Q=\operatorname{diag}(\lambda_1,\ldots,\lambda_n)\) with $\lambda(Q)\in\Gamma_k$. 
At such a point, $D\Phi_k(Q)$ is diagonal, with diagonal entries
\[ 
\frac{\partial}{\partial \lambda_i}\left(-\frac{\sigma_{k-1}}{\sigma_k}\right) = \frac{\sigma_{k-1}\sigma_{k-1,i}-\sigma_k\sigma_{k-2,i}} {\sigma_k^2}.
\]
Here $\sigma_{m,i}$ denotes $\sigma_m$ with $\lambda_i$ omitted, and is understood to be $0$ for $m<0$. The numerator is positive by the standard properties of elementary symmetric polynomials on $\Gamma_k$; see \cite{LT}. Hence $D\Phi_k(Q)>0$. By the real--complex correspondence, $\widetilde\Phi_k$ is elliptic. 

It remains to verify the coefficient control condition. Summing the above entries gives 
\[
\sum_{i=1}^n\frac{\partial}{\partial \lambda_i}\left(-\frac{\sigma_{k-1}}{\sigma_k}\right) = \frac{(n-k+1)\sigma_{k-1}^2-(n-k+2)\sigma_k\sigma_{k-2}} {\sigma_k^2}.
\]
Since $\lambda(Q)\in\Gamma_k$, we have $\sigma_k\sigma_{k-2}\geq 0$. Therefore,
\[
\sum_{i=1}^n\frac{\partial}{\partial \lambda_i}\left(-\frac{\sigma_{k-1}}{\sigma_k}\right) \leq (n-k+1)\left(\frac{\sigma_{k-1}}{\sigma_k}\right)^2.
\]
In terms of $\Phi_k$, the estimate reads
\[ 
\tr\bigl(D\Phi_k(Q)\bigr)\leq (n-k+1)\Phi_k(Q)^2.
\] 
By unitary invariance, the same estimate holds for every Hermitian matrix $Q$ with $\lambda(Q)\in\Gamma_k$. Thus, for any $M\in\sym(2n)$ with $\lambda(\mathcal H(M))\in\Gamma_k$, the trace relation \eqref{tracerelation} gives
\[ 
\tr\bigl(D\widetilde\Phi_k(M)\bigr)\leq\frac{n-k+1}{2}\widetilde\Phi_k(M)^2.
\] 
Thus condition \textup{(iii)} holds with \(G(s)=ns^2\) on $(-\infty,0]$. This verifies the hypotheses of Theorem \ref{Gthm} for $\widetilde\Phi_k$, thereby proving Theorem \ref{CIHQthm}.

\subsection{Deformed Hermitian--Yang--Mills phase operator}
We next verify the hypotheses for this operator. For
$Q\in\herm(n)$, set
\[
\Psi(Q):=\sum_{i=1}^n \arctan \lambda_i(Q),
\]
and let $\widetilde\Psi$ be the induced real operator.

Since
\(
\Psi(Q)<n\pi/2,
\)
condition \textup{(ii)} holds with $c=n\pi/2$. To compute the linearized coefficients, we may again work at a diagonal matrix
$Q=\operatorname{diag}(\lambda_1,\ldots,\lambda_n)$. At such a point,
$D\Psi(Q)$ is diagonal, with diagonal entries
\[
\frac{\partial \Psi}{\partial \lambda_i}=\frac{1}{1+\lambda_i^2}>0.
\]
As in the preceding subsection, this implies the ellipticity of $\widetilde\Psi$. Moreover, using \eqref{tracerelation}, we obtain, for every $M\in\sym(2n)$,
\[
\tr\bigl(D\widetilde\Psi(M)\bigr)
=
\frac{1}{2}\sum_{i=1}^n\frac{1}{1+\lambda_i(\mathcal H(M))^2}\leq\frac{n}{2}.
\]
Therefore condition \textup{(iii)} holds with $G\equiv n/2$ on
$(-\infty,n\pi/2]$. This verifies the hypotheses of Theorem \ref{Gthm}
for $\widetilde\Psi$.

\subsection{Further examples}
Theorem \ref{Gthm} also applies to several other inverse-type operators defined in terms of the eigenvalues of either the real Hessian or the complex Hessian. For example, these include
\[
-\sum_{i=1}^n \lambda_i^{-\alpha},
\quad \alpha>0,
\]
on the positive cone $\Gamma_n$, and
\[
-\left(\frac{\sigma_{k-1}}{\sigma_k}\right)^\alpha,
\quad \alpha>0,
\]
on $\Gamma_k$. Another class is given by inverse $p$-trace-type operators, such as
\[
-\sum_{|I|=p}\frac{1}{\sum_{i\in I}\lambda_i},
\]
on the cone where all $p$-fold sums $\sum_{i\in I}\lambda_i$ are positive. The verification is straightforward and is omitted.

\section{Proof of Theorem \ref{GCIHQthm}}\label{proofofginverseqthm}
\begin{proof}
Following the notation of Subsection \ref{IHQP}, define the induced
real operator
\begin{equation*}
\widetilde\Phi_{k,l}(M)
:=-\frac{\sigma_l(\mathcal H(M))}
{\sigma_k(\mathcal H(M))}.
\end{equation*}
The equation can then be written as
\begin{equation*}
\widetilde\Phi_{k,l}(D^2u)=\frac12 x\cdot Du-u.
\end{equation*}
As in the proof of Theorem \ref{Gthm}, the strategy is to show that
the self-similar term
\begin{equation*}
\phi=\frac12 x\cdot Du-u
\end{equation*}
is constant. 

Write $a^{ij}(x)$ for the $(i,j)$-entry of $D\widetilde\Phi_{k,l}(D^2u(x))$ and define the linearized operator
\begin{equation*}
\mathcal L
:=a^{ij}\partial_{ij}
-\frac12 x\cdot D.
\end{equation*}
Differentiating the equation as in Subsection \ref{subsecEquation for the self-similar term} gives
\begin{equation}\label{eq:GCIHQ-Lphi}
\mathcal L\phi=0.
\end{equation}
It remains to show that $\phi$ attains its global maximum. The strong maximum principle will then imply that $\phi$ is constant. 

As $u$ is $k$-plurisubharmonic (hence subharmonic) on $\mathbb C^n$ and proper, 
\begin{equation*}
u(x)\to +\infty
\qquad\text{as }|x|\to\infty.
\end{equation*} 
We use $u$ itself as a penalizing function and first compute $\mathcal Lu$.

Since $\widetilde\Phi_{k,l}$ is homogeneous of degree $-(k-l)$, Euler's identity gives
\begin{equation*}
a^{ij}u_{ij}=-(k-l)\widetilde\Phi_{k,l}(D^2u)=-(k-l)\phi.
\end{equation*}
Using
\begin{equation*}
\frac12 x\cdot Du=u+\phi,
\end{equation*}
we obtain
\begin{equation}\label{eq:GCIHQ-Lu}
\mathcal Lu=-(k-l)\phi-\frac12 x\cdot Du=-(k-l+1)\phi-u.
\end{equation}

For $\varepsilon>0$, define
\begin{equation*}
w_\varepsilon
:=
\phi-\varepsilon u.
\end{equation*}
Note that
\begin{equation*}
\phi=-\frac{\sigma_l}{\sigma_k}(\partial\bar\partial u)<0.
\end{equation*}
Then
\begin{equation*}
w_\varepsilon(x)\to-\infty
\qquad\text{as }|x|\to\infty.
\end{equation*}
Thus $w_\varepsilon$ attains its global maximum at some point
$p_\varepsilon\in\mathbb C^n$. 

Since $p_\varepsilon$ is a maximizer, $Dw_\varepsilon(p_\varepsilon)=0$, and ellipticity gives
\begin{equation*}
0\geq\mathcal Lw_\varepsilon(p_\varepsilon).
\end{equation*}
By \eqref{eq:GCIHQ-Lphi} and \eqref{eq:GCIHQ-Lu},
\begin{equation*}
0\geq
\mathcal Lw_\varepsilon(p_\varepsilon)=
\varepsilon
\bigl(
u(p_\varepsilon)
+(k-l+1)\phi(p_\varepsilon)
\bigr).
\end{equation*}
Hence
\begin{equation}\label{eq:GCIHQ-maxpoint}
u(p_\varepsilon)
+(k-l+1)\phi(p_\varepsilon)
\leq 0.
\end{equation}

By the maximality of $p_\varepsilon$ and the identity $\phi(0)=-u(0)$,
\begin{equation*}
\phi(p_\varepsilon)
-\varepsilon u(p_\varepsilon)\geq\phi(0)-\varepsilon u(0)=-(1+\varepsilon)u(0).
\end{equation*}
Combining this inequality with \eqref{eq:GCIHQ-maxpoint}, we obtain
\begin{equation*}
\bigl(1+(k-l+1)\varepsilon\bigr)u(p_\varepsilon)\leq(k-l+1)(1+\varepsilon)u(0).
\end{equation*}
It follows that the points $p_\varepsilon$, for
$0<\varepsilon\leq1$, lie in a fixed compact sublevel set of $u$.

Choose a sequence $\varepsilon_j\to 0$. After passing to a subsequence,
we may assume that
\begin{equation*}
p_{\varepsilon_j}\to p_\infty
\end{equation*}
for some $p_\infty\in\mathbb C^n$. For every $x\in\mathbb C^n$, the
maximality of $p_{\varepsilon_j}$ gives
\begin{equation*}
\phi(x)-\varepsilon_j u(x)
\leq
\phi(p_{\varepsilon_j})
-\varepsilon_j u(p_{\varepsilon_j}).
\end{equation*}
Letting $j\to\infty$, we obtain
\begin{equation*}
\phi(x)\leq\phi(p_\infty).
\end{equation*}
Thus $\phi$ attains its global maximum at $p_\infty$. The remainder follows exactly as in the proof of Theorem \ref{Gthm} and is omitted.
\end{proof}
\begin{rem}
The proofs of Theorems \ref{Gthm} and \ref{GCIHQthm} are based on essentially the same underlying idea. In fact, the penalization argument used in the proof of Theorem \ref{GCIHQthm} also yields a shorter proof of Theorem \ref{Gthm}, with $|x|^2$ in place of $u$. On the other hand, the localized barrier argument used for Theorem \ref{Gthm} yields a finite-scale quantitative estimate and may be useful in related problems. Thus the two arguments highlight different aspects of the same method.
\end{rem}

\section*{Statements and Declarations}

\noindent\textbf{Funding.}
The second author was partially supported by the National Natural Science Foundation of China (Grant Nos. 12471054 and 12001292) and the Fundamental Research Funds for the Central Universities, Nankai University (Grant No. 050-63263075).
\medskip

\noindent\textbf{Competing Interests.}
The authors declare that they have no competing interests.

\medskip

\noindent\textbf{Data Availability.}
Data sharing is not applicable to this article, as no datasets were generated or analyzed during the current study.

\medskip

\noindent\textbf{AI usage statement.}
All mathematical ideas, arguments, and results presented in this manuscript were developed independently by the authors, {\bf without the use of generative AI}. The main results were presented by the second author and included in the first author’s master’s thesis, completed in early 2026 for the degree awarded that year. AI tools were used only to polish the language and improve readability. The authors take full responsibility for the content of the manuscript.

\end{document}